\documentclass[a4paper,11pt]{amsart}
\usepackage[T1]{fontenc}
\usepackage[latin1]{inputenc}
\usepackage{enumerate}
\usepackage{amsthm}
\def\car{{\mathbf 1}}

\def\/{\, | \,}

\def\ee{\e}

\def\d{\text{ d}}

\def\Gn{{\Gamma_\Lambda^{(n)}}}%%% Configurations de cardinal n sur Lambda

\def\ee#1{\varepsilon_{_{#1}}}
\def\<{\langle}
\def\>{\rangle}

\def\P{\text{P}}

\def\B{{\mathcal B}}
\def\T{{\mathcal T}}

\def\SS{{\mathfrak S}}

\def\Lc{{\Lambda^c}}
\def\Lnc{{\Lambda_n^c}}

\def\({\Bigl(}
\def\){\Bigr)}
\def\R{{\mathbf R}}
\def\d{\, \text{d}}

\def\hat#1{\ensuremath{#1^s}}
\newcommand{\esp}[1]{{\mathbf E}\left[{#1}\right]}

\newcommand{\Id}{\operatorname{Id}}

\newcommand{\supp}{\operatorname{supp}}

\newcommand{\RE}{{\mathbf R}}
\newcommand{\N}{{\mathbf N}}

\newtheorem{theorem}{Theorem}[section]

\newtheorem{corollary}{Corollary}[section]
\newtheorem{lemma}{Lemma}[section]
\newtheorem{defn}{Definition}
\newtheorem{hyp}{Hypothesis}

\newtheorem{rem}{Remark}[section]

\makeatletter
\let\new@ifnextchar\@ifnextchar
\let\new@ifnch\@ifnch
\makeatother
\begin{document}
\title{Wasserstein distance on configuration space}
\author{L. Decreusefond}
\address{L. Decreusefond \\ GET/ENST - UMR CNRS 5141\\
46, rue Barrault\\ 75634 Paris cedex 13, FRANCE}
 \email{Laurent.Decreusefond@enst.fr}
\begin{abstract}
We investigate here the optimal transportation problem on
configuration space for the quadratic cost. It is shown that, as
usual, provided that the corresponding Wasserstein is finite, there
exists one unique optimal measure and that this measure is supported
by the graph of the derivative (in the sense of the Malliavin
calculus) of a ``concave'' (in a sense to be defined below)
function. For finite point processes, we give a necessary and
sufficient condition for the Wasserstein distance to be finite.
\end{abstract}    
\maketitle

%%%%%%%%%%%%%%%%%%%%%%%%%%%%%%%%%%%%%%%%%%%%%%%%%%%%%%%%%%%%%%%%%%%%%%%%%%%%%%%%%%%%%%%%%%%%%%%%
\section{Introduction}
The optimal transportation problem stems back to the eighteenth
century when G. Monge addressed the optimal way to move earth
particles from one location to another. It is only in the forties of
the last century that Kantorovitch gave this problem its modern form
and a complete solution. According to the Kantorovitch approach, the
optimal transportation problem or Monge-Kantorovitch problem (MKP for
short) reads as follows: given two probability measures $\mu$ and
$\nu$ on a Polish space $X$ and a cost function $c$ on $X\times X$,
does there exist a probability measure $\gamma$ on $X\times X$ which
minimizes $\int c \d \beta$ among all probability measures $\beta$ on
$X\times X$ with first (respectively second) marginal $\mu$
(respectively $\nu$)~? One can furthermore ask whether the optimal
measure is unique and which properties it has. So far, the mainly
investigated situations suppose that  $X=\R^n$ or a finite-dimensional
manifold with a cost function which is $c(x,y)=h(|x-y|)$ where $h$ is
a convex (or concave) function on $\R$.

Varying cost functions and underlying spaces yields to numerous
interesting inequalities often with optimal constants (see
\cite{MR1964483} and references therein) or to new insights on known
theorems such as Strassen Theorem about stochastic ordering (see
\cite{MR99k:28006,MR99k:28007}). Moreover, in a large part of
investigated cases, the optimal measure is unique and is supported by
the graph of a function $T$, i.e., $\gamma=(\Id \otimes
T)^*\mu$. This map $T$ gives raise to a coupling, said
optimal, between the measure $\mu$ and $\nu$: If $A$ is a
r.v. distributed according to $\mu$ then $T(A)$ is distributed
according to $\nu$ and this construction of the two distributions on
the same probability space is optimal in the sense that it minimizes
$E[c(A,B)]$ among all the r.v. $B$ distributed according to
$\nu$. Optimal coupling is a well known tool to obtain inequalities
between random variables (see \cite{MR1741181}) but to the best of our
knowledge, the optimal coupling has to be explicitly built to obtain
these bounds. Meanwhile, optimal transportation theory in the solved
cases, indicates that the optimal coupling (or transportation) map can
be written as the graph of a ``concave'' function, so independently of
the precise description of the map, one can obtain interesting
inequalities just knowing this property of the optimal map.

Our goal is here to develop the theory of optimal transportation for
point processes (or configuration spaces) with the objective to obtain
a machinery yielding inequalities similar to those obtained by
the Stein's method
\cite{MR93g:60043,MR1920275,MR2002k:60029,MR2001e:60040}.

The first step is to define a cost between configurations. Several
possibilities can be envisioned, we chose here a cost with a strong
physical interpretation: given two configurations (or sets of points),
$(x_1, \cdots, \ x_m)$ and $(y_1,\cdots,\ y_m)$, the cost is roughly
defined as $\inf_{\sigma \in \SS_m}\sum_j |x_j-y_{\sigma(j)}|^2$ where
$\SS_m$ is the group of permutations over $\{1,\cdots,\, m\}$ (see
\ref{sec:wass-dist} for the precise definition).  The point is then to
determine the cost to go from a configuration with $m$ points to a
configuration with $m^\prime$ points, when $m\neq m^\prime$. In order
to keep a physical meaning to the definition of our cost, it seems
sensitive to impose an infinite value to something which is
impossible. The negative consequence of this choice is that severe
constraints are imposed (see Theorem \ref{thm:compact-case}) on two finite point
processes for their Wasserstein distance to be finite. On the other
hand, the positive consequence is that the optimal measure has a well
defined structure. These constraints disappear when we deal with
locally finite but not finite point processes and we still have a
rigid structure for the optimal measure. It turns out that proving
 here the uniqueness and describing the form of the optimal measure
 is highly similar to the proof of the same properties for the Wiener
 space (see \cite{MR2036490}).

 This paper is organized as follows : we  describe the
 Monge-Kantorovitch Problem in its general settings for a generic cost
 function on a product of two abstract Polish spaces, since we will
 need to instantiate these general results to different particular
 situations in the subsequent sections. Section \ref{sec:wass-dist} is
 devoted to general properties of the Wasserstein distance on
 configuration spaces irrespective to the properties of finiteness of
 the considered point processes. In Section
 \ref{sec:finite-point-proc}, we work under the assumption that only a
 finite number of atoms are random in the $\mu$-configurations and we
slightly  modify  our cost function to simultaneously solve optimal
 transportation between finite point processes and pave the way to the
 analysis for locally finite point processes. This latter analysis is
 done in Section \ref{sec:locally-finite-point}.
%%%%%%%%%%%%%%%%%%%%%%%%%%%%%%%%%%%%%%%%%%%%%%%%%%%%%%%%%%%%%%%%%%%%%%%%%%%%%%%%%%%%%%%%%%%%%%%% 
\section{Preliminaries}
For $X$ and $Y$ two Polish spaces, for $\mu$ (respectively $\nu$) a
probability measure on $X$ (respectively $Y$), $\Sigma(\mu, \, \nu)$
is the set of probability measures on $X\times Y$ whose first marginal
is $\mu$ and second marginal is $\nu$. We also need to consider a
lower semi continuous function $c$ from $X\times Y$ to $\R^+$. The
Monge-Kantorovitch problem associated to $\mu$, $\nu$ and $c$, denoted by MKP($\mu$, $\nu$, $c$) for short, consists in finding 
\begin{equation}
\label{eq:3}
  \inf_{\gamma\in \Sigma(\mu, \, \nu)}\int_{X\times Y}c(x,y)\d \gamma(x,\, y).
\end{equation}
More precisely, since $X$ and $Y$ are Polish and $c$ is l.s.c., it is
known from the general theory of optimal transportation, that there
exists an  optimal measure $\gamma\in \Sigma(\mu,\, \nu)$ and that the minimum coincides with
\begin{equation*}
  \sup_{(F,\, G)\in \Phi_c}(\int_X F\d \mu+\int_Y G\d \nu),
\end{equation*}
where $(F,\, G)$ belongs whenever $F\in L^1(\d \mu)$, $G\in L^1(\d
\nu)$ and $F(x)+G(y)\le c(x,\, y)$.  We will denote by $\T_c(\mu, \,
\nu)$ the value of the infimum in \eqref{eq:3}.  Solving the
Monge-Kantorovitch problem on Polish spaces, is then essentially
proving the finiteness of (\ref{eq:3}) and the uniqueness of the
optimal measure. For $X=Y=\R^k$ and $c$ taken to be the square euclidean
distance, we have the second moment condition: for $\T_c(\mu,\, \nu)$
to be finite it is sufficient that $\int \|x\|^2\d \mu$ and $\int
\|x\|^2\d \nu$ are finite, where $\|x\|$ is the euclidean norm (see
\cite{MR1964483}). For further reference, we denote by $\T_e$ (where
$e$ stands for euclidean) this
distance between probability measures on $\R^k$. If $\mu$ is the Gaussian measure on $\R^k$
and $\nu$ is absolutely continuous with respect to $\mu$ with
Radon-Nikodym density $L$, it is sufficient that $L$ has a finite
entropy, i.e., the $\mu$-expectation of $L\ln L$ is finite, for
$\T_e(\mu,\, L\mu)$ to be finite. This criterion extends to the infinite dimensional setting
where $X=Y$ is a Wiener space and $c(x,y)=2^{-1}\|x-y\|_H^2$, where
$H$ is the associated Cameron-Martin space (see \cite{MR2036490}). In full
generality, we know from \cite{MR99k:28006} that if there exist $F\in
L^1(\d\mu)$ and $G\in L^1(\d\nu)$ such that $c(x,y)\le F(x)+G(y)$,
then $\T_c(\mu, \, \nu)$ is finite.

Once the finiteness of $\T_c$ is ensured, it remains to know whether
the optimal measure is unique. For, it is essential to see that a
measure $\gamma \in \Sigma(\mu, \, \nu)$, is optimal if and only
if its support  is $c$-cyclically monotone (see \cite{MR1699061,MR1395577}) : for any $((x_i,\, y_i),\, i=1,\cdots,m)\in (\supp \gamma)^m$,  we have
\begin{equation*}
  \sum_{i=1}^m c(x_i,\, y_i)\le \sum_{i=1}^m c(x_i,\, y_{\sigma(i)}),
\end{equation*}
for any $\sigma \in \SS_m$, the group of permutations over $\{1, \cdots,\, m\}$.
Moreover, the support of any optimal measure is included in the $c$-super-gradient of a $c$-concave function: 
For $F\, :\, X \to \bar{\RE}=\R\cup \{+\infty\},$
its $c$-super-gradient, $ \partial^c F,$  is the subset of $\Gamma_X\times\Gamma_Y$ of
$(x,\, y)$ such that $c(x,\, y )< \infty$ and 
\begin{displaymath}
  \ F(x)-F(z)\ge
  c(x,\, y)-c(z,\, y), \text{ for any } z \text{ such that } F(z)<+\infty.
\end{displaymath}
The section at $x,$ $ \partial^cF(x)$, is the set $\{y\in Y, \, (x,\,
y)\in \partial^cF\}$. A function $F\, :\, X \to \RE$ is called
c-concave if there exist a set index $I$, $(y_i, \, i\in I)$ a family
of elements of $Y$ and $(a_i, \, i\in I)$ a family of real numbers
such that
\begin{equation*}
%  \label{eq:12}
  F(x)=\inf_{i\in I}(c(x,y_i)+a_i).
\end{equation*}
If we prove that the $c$-super gradient of a $c$-concave function is
single valued, we are done, i.e., we have proved the uniqueness of the
optimal measure. Indeed, if $\partial^cF(x)$ is reduced to a singleton
for $\mu$-a.s. any $x\in X$, this means that $\partial^cF(x)$ is
closed and the selection theorem then induces that there exists a
measurable map $T$ such that $(x,y)$ belongs to $\partial^cF$ if and
only if $y=T(x)$. The uniqueness follows then from the following lemma which we borrow from
\cite{MR2036490}.
\begin{theorem}[See \protect\cite{MR2036490}]
  \label{thm:single_valued_map_implies_uniqueness}
Let $X$ and $Y$ be two Polish spaces and $c$ be a
lower-semi-continuous function from $X\times Y$ to
$\RE^+\cup\{+\infty\}$. Consider the Monge-Kantorovitch problem
associated to the marginals $\mu$ on $X$, $\nu$ on  $Y$ and~$c$. Assume that for any optimal measure $\gamma$,
there exists a measurable map $T_\gamma$ such that $\text{supp }\gamma\subset
\{(x,\, T_\gamma(x)),\, x\in \text{supp }\mu\}$. Then, there exist a
unique optimal measure $\gamma$ and a unique measurable map $T$ such that $\gamma=(\Id \otimes T)^* \mu$.
\end{theorem}
\begin{proof}
For any probability measure $\gamma$ on $X\times Y$, we denote by $J(\gamma)$  the integral of $c$ with respect to $\gamma$: 
  \begin{equation*}
    J(\gamma)=\int_{X\times Y}c(x,\, y)\d \gamma(x,\, y).
  \end{equation*}
  Assume that $\gamma_1$ and $\gamma_2$ are two different optimal
  measures. 
Since $J$ is linear with respect to $\gamma$,
  $\gamma_0=(\gamma_1+\gamma_2)/2$ is also optimal. We denote by $T_0$
  a map whose graph contains the support of $\gamma_0$. Furthermore, for
  $i=1,2$,   $\gamma_i$ is absolutely continuous with respect to
  $\gamma_0$. We denote by $L_i$ the Radon-Nikodym derivative of
  $\gamma_i$ with respect to $\gamma_0$. For any $f\in {\mathcal
    C}_b(X)$, we have
  \begin{align*}
    \int_X f(x)\d\mu(x) & = \int_{X\times Y} f(x)\d\gamma_i(x,\,y)\\
& = \int_{X\times Y} f(x)L_i(x,y)\d\gamma_0(x,\, y)\\
&=\int_{X}f(x)L_i(x,\, T_0(x))\d\mu(x).
  \end{align*}
Therefore, we must have $L_i(x,\, T_0(x))=1$ $\mu$-a.s. or in other
words, $L_i=1$ $\gamma_0$-a.s. for $i=1,2.$ This means that
$\gamma_1=\gamma_2$ and then the uniqueness of the optimal measure for
MKP($\mu,\, \nu, \ c$).

Assume now that there exist two maps $T_1$ and $T_2$ such that
$\gamma=(\Id \otimes T_i)^*\mu.$ This implies that for any $f\in
{\mathcal C_b}(X)$, $f\circ T_1=f\circ T_2$ $\mu$-a.s. hence that
$T_1=T_2$ $\mu$-almost surely. 
\end{proof}
The simplest way to prove that the $c$ super-gradient of a $c$-concave
function is single-valued is to show that a $c$-concave is
``differentiable'' in some sense. That is why, we need to introduce a notion of gradient on configuration space.
The notations are mainly those of  \cite{MR99d:58179}. 
Let $\Gamma_X$ be the  configuration space over a Polish space $X,$ i.e., 
\begin{displaymath}
  \Gamma_X=\{\eta \subset X; \ \eta\cap K \text{ is a finite set
    for every compact } K\subset X \}.
\end{displaymath}
 We identify $\eta\in\Gamma_X$ and the positive Radon measure
$\sum_{x\in\eta} \varepsilon_x.$ Throughout this paper, $\Gamma_X$ is endowed with the vague
topology,  i.e., the weakest topology such that for all $f\in {\mathcal C}_0$ (continuous with compact support on
$X$), the maps 
\begin{displaymath}
  \eta \mapsto \int f\d \eta=\sum_{x\in\eta} f(x)
\end{displaymath}
are continuous. When $f$ is the indicator function of a subset $B,$ we
will use  the shorter  notation $\eta( B)$ to denote the integral of $\car_B$ with respect to $\eta$.  We denote by ${\mathcal B}(\Gamma_X)$ the
corresponding Borel $\sigma$-algebra.

The intensity measure of a probability measure $\mu$ on $\Gamma_X,$ is denoted by $E_\mu\eta$ and defined by
  \begin{math}
    (E_\mu\eta)(B)=E_\mu[\eta( B)],
  \end{math}
for any $B\in {\mathcal B}(\Gamma_X).$ 
We assume henceforth that $E_\mu\eta$ is a positive Radon
 measure on ${\mathcal B}(\Gamma_X).$ 

In what follows, we will take $X=\R^k$ for some $k\ge 1$.
% Let $\Diff$ be the group of diffeomorphisms of $X$ and denote by
% $\DiffO$ the subgroup of all diffeomorphisms with compact support,
% i.e., equal to the identity outside a compact set.
Let $V(X)$ be the set of ${\mathcal C}^\infty$ vector fields on
$X$ and $V_0(X)\subset V(X),$ the subset consisting of all
vector fields with compact support.
For $v\in V_0(X),$ for any $x\in X,$ the curve 
\begin{displaymath}
  t\mapsto {\mathcal V}_t^v(x)\in X
\end{displaymath}
is defined as the solution of the following Cauchy problem
\begin{equation}
  \label{eq:4}
  \begin{cases}
    \dfrac{d}{dt}{\mathcal V}_t^v(x)&=v({\mathcal V}_t^v(x)),\\
{\mathcal V}_0^v(x)&=x.
  \end{cases}
\end{equation}
The associated flow $({\mathcal V}_t^v,\, t\in \R)$ induces a curve
$({\mathcal V}_t^v)^*\eta=\eta\circ ({\mathcal V}_t^v)^{-1}$, $t\in \R$, on
$\Gamma_X$: If $\eta=\sum_{x\in \eta}\varepsilon_x$ then  $({\mathcal
   V}_t^v)^*\eta=\sum_{x\in \eta}\varepsilon_{{\mathcal V}_t^v(x)}.$

 \begin{hyp}
\label{hyp:hyp}
Throughout this paper, we assume that $\mu$ (or $\nu$) is a Borel probability measure on $\Gamma_X$ such that the following conditions hold.
\begin{enumerate}[i) ]
\item \label{item:1}$\eta(\{x\})\in \{0,\, 1\}$ for all $x\in X$ and $\mu$-a.s. $\eta$.
\item  \label{item:2} Either $\mu(\eta: \, \eta(X)<+\infty)=1$ or $\mu(\eta: \, \eta(X) = +\infty)=1$.
\item \label{item:3} For all $v \in V_0(X)$ and $t\in \R$, $\mu$ is quasi-invariant with respect to the flow $({\mathcal V}_t^v)^*$ of $\Gamma_X$, i.e., $\mu\circ \left(({\mathcal V}_t^v)^*\right)^{-1}$ is equivalent to $\mu$.
\end{enumerate}   
 \end{hyp}
We are then in position to define the notion of differentiability on $\Gamma_X$.
A measurable function $F\, : \, \Gamma_X \to \R$ is said to be differentiable if for any $v\in V_0(X)$, the following limit exists:
\begin{equation*}
 \lim_{t\to 0} t^{-1}\left( F(({\mathcal V}_t^v)^*\eta)-F(\eta)\right).
\end{equation*}
We then denote $\nabla^\Gamma_vF(\eta)$ the preceding limit and by
$\nabla_x^\Gamma F(\eta)$ the corresponding gradient (see \cite{MR99d:58179}) which is
defined by the identity:
\begin{equation*}
  \int \nabla_x^\Gamma F(\eta).v(x)\d \eta(x)=\nabla^\Gamma_vF(\eta),
\end{equation*}
for all $v\in V_0(X)$.

%%%%%%%%%%%%%%%%%%%%%%%%%%%%%%%%%%%%%%%%%%%%%%%%%%%%%%%%%%%%%%%%%%%%%%%%%%%%%%%%%%%%%%%%%%%%%%%%

\section{Wasserstein distance}
\label{sec:wass-dist}

We consider on $X=\R^k$ the cost function as $d(x,y)=2^{-1}\|x-y\|^2$ where $\|x\|$ denote the euclidean norm of $x\in X$
and we  define a cost between configurations (see also
\cite{MR1920275,MR2001e:60040,MR2001a:60058}) as the 'lifting' of $d$
on $\Gamma_X$:
\begin{equation*}
%\label{eq:definition_du_cout}
  c(\eta_1,\eta_2)=\inf \left\{ \int d(x,y) d\beta(x,y),\
    \beta\in \Gamma_{\eta_1,\eta_2}\right\},
\end{equation*}
where $\Gamma_{\eta_1,\eta_2}$ denotes the set of
$\beta\in\Gamma_{X\times X}$ having marginals $\eta_1$ and $\eta_2.$
According to \cite{MR1730565}, $c$ is lower semi continuous on
$\Gamma_{X}\times\Gamma_X$. 
We can then set the Monge-Kantorovitch
problem for configuration spaces.
\begin{defn}
  Let $\mu$ and $\nu$ be two probability measures on
  $(\Gamma_X,{\mathcal B}(\Gamma_X))$. We say that a probability
  $\gamma$ on $(\Gamma_X\times \Gamma_X,{\mathcal B}(\Gamma_X\times
  \Gamma_X))$ is a solution of the Monge-Kantorovitch Problem
  associated to the couple $(\mu,\nu)$ and to the cost $c$ (MKP($\mu,\nu,c$) for
  short) if the first marginal of $\gamma$ is $\mu,$ the second one
  is $\nu$ and if
  \begin{equation*}
    \begin{split}
      J(\gamma)& =\int c(\eta,\zeta) \ d\gamma(\eta,\, \zeta)\\
      & =\inf\left\{\int c(\eta,\zeta) \ d\beta(\eta,\, \zeta)\, :\,
        \beta\in \Sigma(\mu,\nu)
\right\}
    \end{split}
  \end{equation*}
\end{defn}
Since $\Gamma_X$ is  Polish, this infimum is attained and is equal to 
\begin{displaymath}
  \sup\{\int F(\eta) \d\mu(\eta) +\int G(\zeta) \d\nu(\zeta)\,
    :\, (F,G)\in\Phi_c\},
\end{displaymath}
where $\Phi_c$ is the set of pairs of measurable, real-valued
functions $F$ and $G$ such that $F$ (resp. $G$) belongs to
$L^1(\d\mu)$ (resp. $L^1(\d\nu)$) and $F(\eta)+G(\zeta)\le c(\eta,\zeta).$ The Wasserstein distance between $\mu$ and $\nu$ is the square root of $\T_c(\mu, \, \nu)$.

Since the cost $c$ is infinite whenever the two
configurations do not have the same mass, we have the following
theorem.
\begin{theorem}
\label{thm:egalite_loi_nb_atomes}
Let $\mu$ and $\nu$ be two probability measures on the
configuration space $\Gamma_X$. If the Monge-Kantorovitch cost, with respect
to $c$, is finite then
\begin{displaymath}
\mu(\eta(X)=n)=\nu(\eta(X)=n) \text{ for any } n\in \N\cup \{+\infty\}.  
\end{displaymath}
\end{theorem}
\begin{proof}
    There exists
  at least one  measure $\gamma$ such that 
  \begin{multline*}
    \infty > \int c(\eta,\omega)\d\gamma(\eta,\omega) \\ =
    \int_{\eta(X)=\omega(X)}c(\eta,\omega)\d\gamma(\eta,\omega)+\int_{\eta(X)\neq \omega(X)}c(\eta,\omega)\d\gamma(\eta,\omega).  \end{multline*}
This implies that $\gamma(\eta(X)\neq \omega(X))=0$. It follows that 
\begin{align*}
  \nu(\omega(X)=n)&=\gamma(\omega(X)=n)\\
 &= \gamma(\omega(X)=n;\, \eta(X)=n) +
  \gamma(\omega(X)=n;\, \eta(X)\neq n)\\
&= \gamma(\omega(X)=n;\, \eta(X)=n),
\end{align*}
for any  $n\in \N\cup \{+\infty\}$. By the very same reasoning, it also holds that 
$\mu(\eta(X)=n)=\gamma(\omega(X)=n;\, \eta(X)=n)$ for any $n$ and
thus that $\mu(\eta(X)=n)=\nu(\eta(X)=n)$ for any $n\in \N\cup \{+\infty\}$.
\end{proof}
%%%%%%%%%%%%%%%%%%%%%%%%%%%%%%%%%%%%%%%%%%%%%%%%%%%%%%%%%%%%%%%%%%%%%%%%%%%%%%%%%%%%%%%%%%%%%%%%%%%%%%%%%%%%%%%%%%%%%%%%%%%%%%
\section{Finite point processes}
\label{sec:finite-point-proc}
Consider $\Lambda$ a compact set of $X$ and let $\zeta$ be fixed in
$\Gamma_{\Lambda^c}$. We define $c_\zeta$ as:
\begin{equation*}
  \begin{array}{llll}
    c_\zeta\, :\, & \Gamma_\Lambda \times \Gamma_X &\to
  &\RE^+\cup\{+\infty\}\\
& (\eta,\, \omega) &\mapsto &c(\eta+\zeta,\, \omega).
  \end{array}
\end{equation*}
Let $\mu$ be a probability measure on $\Gamma_\Lambda$ and $\nu$ a probability measure on $\Gamma_X$, we denote by $\T_{c_\zeta}(\mu,\, \nu)$ the $c_\zeta$-Wasserstein distance between $\mu$ and $\nu$:
\begin{equation*}
  \T_{c_\zeta}(\mu,\, \nu)=\inf_{\gamma \in \Sigma(\mu,\, \nu)}\int_{\Lambda\times X} c_\zeta(\eta, \, \omega)\d \gamma(\eta,\, \omega). 
\end{equation*}
Since $\Lambda$ is compact and $E_\mu\eta$ is supposed to be a Radon measure, the configurations of $\Gamma_\Lambda$ have $\mu$-a.s. a finite number of atoms. It it then useful to think of $\Gamma_\Lambda$ as the disjoint union of the spaces $\Gamma_\Lambda^{(n)}$ for $n$ running from $0$ to infinity, where
\begin{equation*}
  \Gamma_\Lambda^{(n)}=\{\eta \in \Gamma_\Lambda, \ \eta(\Lambda)=n\}.
\end{equation*}
Then, consider $\tilde{\Lambda}^n=\{(x_1,\cdots,\, x_n)\in \Lambda^n,\, x_i\neq x_j \text{ for } i\neq j\}$, there is a bijection $s_\Lambda^n$ between $\tilde{\Lambda}^n/\SS_n$ and $\Gamma_\Lambda^{(n)}$:
  \begin{eqnarray*}
  s_\Lambda^n :&  \tilde{\Lambda}^n/\SS_n   &\xrightarrow{\hphantom{\longmapsto}} \Gamma_\Lambda^{(n)}\\
 & \{x_1,\cdots,x_n\} &\longmapsto \sum_{i=1}^n \ee{x_i}.
\end{eqnarray*}
The topology of $\tilde{\Lambda}^n/\SS_n$ induced by the usual
topology of $\Lambda^n$ thus defines a locally compact metrizable
Hausdorff topology on $\Gamma_\Lambda^{(n)}$. Since $\Lambda$ is
compact, this topology coincides with the restriction to
$\Gamma_\Lambda^{(n)}$ of the vague topology on $\Gamma_\Lambda$.  We put on
$\Gamma_\Lambda^{(n)}$, $\B(\Gamma_\Lambda^{(n)})$ the associated
Borel $\sigma$-algebra.  For any map $F$ from $\Gamma_\Lambda$
into a measurable space $(Y,{\mathcal Y})$, for any integer $n$, we can
consider, $F_n$, the restriction of  $F_n$ to $\Gn$:
\begin{equation*}
   \begin{array}{llll}
    F_n\, :\, & \Gamma_\Lambda^{(n)}  &\to
  &Y\\
& \eta &\mapsto &F_n(\eta)=F(\eta).
  \end{array}
\end{equation*}
Since $\Gamma_\Lambda^{(n)}$ is closed in $\Gamma_\Lambda$, it is a Polish space and  $F_n$
is measurable from $(\Gn, \, {\mathcal B}(\Gn))$ into $(Y,\, {\mathcal
  Y}).$

 We now identify $\sigma \in \SS_n$ and its
action over $\tilde{\Lambda}^n$, which maps $x=(x_1,\cdots,\, x_n)$ to $\sigma
x=(x_{\sigma(1)},\cdots, \, x_{\sigma(n)})$. Let $\hat{F}_n$ be a measurable function from $\tilde{\Lambda}^n$
into a measurable space $(Y,\,{\mathcal Y})$.  We say that $\hat{F}_n$
is symmetric whenever for any $\sigma \in \SS_n$, $\hat{F}(\sigma
x)=\hat{F}(x)$ for any $x\in \tilde{\Lambda}^n$.  Identify now $\tilde{\Lambda}^n
/\SS_n$ with a subset $\Lambda^\prime$ of $\tilde{\Lambda}^n$, since
$\tilde{\Lambda}^n$ has $n!$ disjoint connected components, the map
\begin{equation*}
  \begin{array}{lccl}
    j^n \, : \, & \tilde{\Lambda}^n & \longrightarrow & \tilde{\Lambda}^n/\SS_n\times \SS_n\\
&&&\\
 & x & \longmapsto & (\bar{x},\, \sigma)=(j^n_1(x), \, j^n_2(x)),
  \end{array}
\end{equation*}
where $\sigma$ is such that $\sigma x=\bar{x}\in \Lambda^\prime$, is
an homeomorphism. Furthermore, $j^n_1$ is a local
diffeomorphism. Hence, any symmetric measurable (respectively continuous
or differentiable) function $\hat{F}_n$ from $\tilde{\Lambda}^n$ into $\R$ can
be identified with a measurable (respectively continuous or
differentiable) function $F_n$ from $\Gamma_\Lambda^{(n)}$ into
$\R$ with $$F_n(\eta)=\hat{F}_n\left(j^{-1}((s_\Lambda^n)^{-1}(\eta), \,
  \sigma)\right)$$ for any $\sigma\in \SS_n$ or equivalently
with $$F_n(\eta)=\frac{1}{n!}\sum_{\sigma\in
  \SS_n}\hat{F}_n\left(j^{-1}((s_\Lambda^n)^{-1}(\eta), \, \sigma)\right).$$
Conversely, any function $F_n$ from $\Gamma_\Lambda^{(n)}$ into
$\R$ gives raise to a symmetric function $\hat{F}_n$ from $\tilde{\Lambda}^n$
into $\R$ by $\hat{F}_n(x)={F}(s_\Lambda^n(j_1(x)))$, with the same
regularity (measurable, continuous or differentiable). Accordingly,
every probability measure $\mu_n$ on $\Gamma_\Lambda^{(n)}$ can be
viewed as a symmetric (i.e., invariant under the action of $\SS_n$)
probability measure $\mu^s_n$ on $\tilde{\Lambda}^n$ and vice-versa.

  Let $\mu$ be a probability measure on $\Gamma_\Lambda$ and consider the
  disintegration of $\mu$ along the map $(\eta \mapsto \eta(\Lambda))$:
  \begin{equation*}
    \mu(B)=\sum_{n\ge 0} \mu(B\/ \eta(\Lambda)=n) \P(\eta( \Lambda)=n).
  \end{equation*}
  We denote by $\mu_n$ the measure $\mu(.\/ \eta(\Lambda)=n)$. The
  measure $\P(\eta( \Lambda)=n)\mu_n$ is the so-called Janossy
  measure of order $n$ (see \cite{MR1950431}).  We say that
  $\mu$ is regular whenever for any $n\ge 1$, $\hat{\mu}_n$, the symmetric measure
  associated to $\mu_n$, is absolutely continuous with respect to the Lebesgue
  measure on $\tilde{\Lambda}^n$. 
\begin{rem}
\label{rem:finite-point-proc}
  Since $X$ is Polish it can embedded as a $G_\delta$ in a 
  compact metric space $X^\prime$. If a probability measure $\nu$ on
  $\Gamma_X$ is such that $\nu(\omega(X)<+\infty)=1$,  we can embed
  $(\Gamma_X, \B(\Gamma_X),\, \nu)$ into $(\Gamma_{X^\prime},
  \B(\Gamma_{X^\prime}), \, \nu_{X^\prime})$, with $\supp
  \nu_{X^\prime}=\supp \nu$ and
  $\nu_{X^\prime}(\omega(X)<+\infty)=1$. Thus, all the previous results
  established on $\Gamma_\Lambda$ are valid on $\Gamma_{X^\prime}$
  hence on $\Gamma_X$. In particular to every probability measure
  $\nu_n$ on $\Gamma_X\subset \Gamma_{X^\prime}$, we can associate, as
  above, a  
  symmetric probability measure, $\hat{\nu}_n$ on $\tilde{X}^n$.
\end{rem}
The next theorem follows from the previous considerations.
\begin{theorem}
\label{thm:equi_differentiabilite}
  Assume that $\mu$ is a regular probability measure on
  $\Gamma_\Lambda$ and let $F$ be  measurable from $\Gamma_\Lambda$
  into $\R$. Then, $F$ is $\mu$-a.s. differentiable, on its domain, if and only if $\hat{F_n}$ is
  $\hat{\mu_n}$-a.s. differentiable, on its domain, for any integer $n$. 
\end{theorem}
 The euclidean symmetric cost on $X^n$, denoted by $c^s_n$, is defined as:
 \begin{equation*}
   c^s_n(x,\, y)=\inf_{\sigma \in \SS_n}\frac{1}{2}\| x-\sigma y\|^2.
 \end{equation*}
It is immediate that 
 \begin{equation}
\label{eq:5}
   c^s_n(x,\, y)=c\left(s_\Lambda^n(j_1^n (x)),\, s_\Lambda^n(j_1^n(y)) \right)
 \end{equation}
 and that 
\begin{equation}\label{eq:6}
  c(\eta,\, \omega)=c^s_n(x,\, y)
\end{equation}
for any $x\in (s_\Lambda^n\circ j_1^n)^{-1}(\{\eta\})$ and any $y\in (s_\Lambda^n\circ j_1^n)^{-1}(\{\omega\})$.
\begin{lemma}
  \label{thm:cs-concave-equivaut-concave-symm}
A function $\hat{F}_n$ from $\tilde{\Lambda}^n$ into $\R$ is $c^s_n$-concave if
and only if $\hat{F}_n-\|x\|^2/2$ is concave in the usual sense and
$\hat{F}_n$ is symmetric.
\end{lemma}
\begin{proof}
By its very definition, a $c^s_n$-concave function  $\hat{F}_n$ is of the form:
  \begin{equation}
\label{eq:euc-concave}
    \hat{F}_n(x)=\inf_{i\in I} (c^s_n(x,\, y_i)+a_i) = \inf_{\substack{i\in I\\
        \sigma \in \SS_n}}(\frac{1}{2}\| x-\sigma y_i\|^2 + a_i),
  \end{equation}
where $y_i$ belongs to $\tilde{\Lambda}^n$ for any $i\in I$.
This clearly implies that $\hat{F}_n$ is symmetric and
euclidean-concave and  euclidean-concavity is known to be equivalent to
the concavity of $(x\mapsto \hat{F}_n(x)-\|x\|^2/2)$  in the usual
sense (see \cite{MR1964483}), hence the result. 

It only remains to prove that $\hat{F}_n$ symmetric and
euclidean-concave can be written as in \eqref{eq:euc-concave}. Since
$\hat{F}_n$ is euclidean concave,  
\begin{displaymath}
  \hat{F}_n(x)=\inf_{i\in I}(\frac{1}{2}\| x-y_i\|^2 + a_i),
\end{displaymath}
for some index set $I$, $(a_i,\, i\in I)$ a family of real numbers and
$(y_i,\, i\in I)$ some elements of $X^n$.
Since $\hat{F}_n$ is symmetric, $\hat{F}_n(x)=\hat{F}_n(\sigma x)=\inf_{\sigma\in \SS_n}\hat{F}_n(\sigma (x))$ thus 
\begin{displaymath}
  \hat{F}_n(x)=\inf_{\substack{i\in I\\ \sigma\in \SS_n}}(\frac{1}{2}\| \sigma
  x-y_i\|^2 + a_i)=\inf_{\substack{i\in I\\ \sigma\in
      \SS_n}}(\frac{1}{2}\|x- \sigma y_i\|^2 + a_i).
\end{displaymath}
% Since $\tilde{\Lambda}^n$ is dense in $\Lambda^n$, each $y_i\in \Lambda^n\backslash \tilde{\Lambda}^n$ can be replaced by a sequence $(y_i^k,\, k\ge 1)$ of elements of $\tilde{\Lambda}^n$, converging to $y_i$, without altering the value of the infimum:
% \begin{equation*}
%   \|x- \sigma y_i\|^2 =\inf_k\|x- \sigma y_i^k\|^2.
% \end{equation*}
The proof is thus complete.
\end{proof}
It follows from the Lebesgue-a.s. differentiability of concave
function that we have:
\begin{corollary}
  \label{thm:differentiabilite-cs-concave}
Let $n\ge 1$ and  $\hat{F}_n$ be a $c^s_n$-concave function. Then, $\hat{F}_n$ is
Lebesgue-a.s. differentiable on its domain. 
\end{corollary}
\begin{corollary}
\label{cor:single-valued}
  Let $n\ge 1$, $\hat{\mu}_n$ an absolutely continuous measure on $\tilde{\Lambda}^n$ and  $\hat{F}_n$  a $c^s_n$-concave function. Then,
  $\partial_{c_n^s}\hat{F}$ is $\hat{\mu}_n$-a.s. single-valued.
\end{corollary}
\begin{proof}
  We already know (see Corollary \ref{thm:differentiabilite-cs-concave}) that
  $\hat{F}_n$ is Lebesgue-a.s. differentiable. From (\ref{eq:6}), it
  is clear that 
  \begin{equation*}
    \partial^{c_n^s}\hat{F}_n(x)=\partial^{c_\emptyset}F_n\bigl(s_\Lambda^n(j^n_1(x))\bigr).
  \end{equation*}
The previous theorem implies that the rightmost set is reduced to a
singleton for $\hat{\mu}_n$-almost-all $x$, hence the result. 
\end{proof}
Remind now that for two configurations $\eta+\zeta$ and $\omega$ at finite
$c$ distance, $\beta_{\eta+\zeta,\, \omega}$ is one measure on
$\Gamma_{X\times X}$ which realizes this distance. 
\begin{defn}
For any $\eta\in \Gamma_X$, for any $\Lambda\subset X$,
$\pi^\Lambda(\eta)=\eta\cap \Lambda.$ For any  map $t$ from $X$ to
$X$, we associate the map $t^\Gamma$ from $\Gamma_X$ to itself,  defined by
\begin{equation*}
  t^\Gamma(\sum_{x\in \eta}\varepsilon_x)=\sum_{x\in \eta}\varepsilon_{t(x)}
  \text{ for any } \eta=\sum_{x\in \eta} \varepsilon_x
\end{equation*}
For any probability measure
$\mu$ on
$\Gamma_X$, $\pi^\Lambda\mu$ is the image measure of $\mu$ by $\pi^\Lambda$.
For any $\eta=(\eta_1,\, \eta_2) \in \Gamma_X\times \Gamma_X,$ we set
$p_i(\eta)=\eta_i$ for $i=1,\, 2$. Accordingly, for any probability
measure $\gamma$ on $\Gamma_X \times \Gamma_X$, $p_i\gamma$ the image
of $\gamma$ by $p_i$. We also introduce  $\pi_i^\Lambda:=\pi^\Lambda\circ
p_i$, thus $\pi_1^\Lambda(\eta,\, \omega)$ is the restriction to
$\Lambda$ of $\eta$.
  For any configuration $\beta$ on $X\times X$, define  $r^\Lambda$ by:
  \begin{equation*}
   \begin{array}{llll}
    r^\Lambda \, :\, & \Gamma_{X\times X} & \to & \Gamma_{X\times X} \\
& \beta & \mapsto & r^\Lambda\beta=\beta \cap (\Lambda\times X).
  \end{array} 
  \end{equation*}
At last,  $r^\Lambda_i$ denotes
$p_i\circ r^\Lambda$.
\end{defn}
The main result of this section is the following.
\begin{theorem}
\label{thm:compact-case}
Let $\mu$ (resp. $\nu$) be a probability  measure on $\Gamma_\Lambda$
(resp. $\Gamma_X$) satisfying Hypothesis \ref{hyp:hyp} and $\zeta\in \Gamma_{\Lambda^c}$. Assume that $\mu$ is regular and that
$\T_{c_\zeta}(\mu,\, \nu)$ is finite. Then, for any optimal measure
$\rho$, there exists a $c_\zeta$-concave function $F$ such that $\supp
\rho \subset \partial^{c_\zeta} F$ and for any $\omega\in
\partial^{c_\zeta} F(\eta)$,   
\begin{equation*}
    r^\Lambda(\beta_{\eta+\zeta,\, \omega})=\sum_{x\in\eta}\varepsilon_{(x,\,
    x-\nabla^\Gamma _x F(\eta))},
\end{equation*}
 for any
$\beta_{\eta+\zeta,\, \omega}$ realizing $c(\eta+\zeta,\, \omega).$
\end{theorem}
\begin{proof}
  $\Gamma_\Lambda$ and $\Gamma_X$ are Polish spaces hence there exists
  at least an optimal measure $\rho$ and a $c_\zeta$-concave function
  $F$ such that $\supp \rho \subset \partial^{c_\zeta} F$. By the
  definition of $c_\zeta$-concavity, for any $\eta\in \Gamma_\Lambda$,
  \begin{align*}
    F(\eta)&=\inf_{i\in I} (c(\eta+\zeta,\ \omega_i) + a_i)\\
&=\inf_{i\in I}\inf_{\varpi_i \subset \omega_i}(c(\eta,\, \varpi_i)+c(\zeta,\varpi_i^c)+a_i)
  \end{align*}
Since $c(\zeta,\varpi_i^c)+a_i$ does not depend on $\eta$, $\hat{F}_n$
is $c^s_n$-concave for any integer $n$. Then Corollary \ref{thm:differentiabilite-cs-concave} implies that $\hat{F}_n$
is Lebesgue-a.s. differentiable, which in turn entails that $\hat{F}_n$
is $\hat{\mu}_n$-a.s. differentiable, since $\hat{\mu}_n$ is absolutely
continuous with respect to the Lebesgue measure. Thus, according to Corollary \ref{thm:differentiabilite-cs-concave} and
Theorem \ref{thm:equi_differentiabilite}, $F$ has $\mu$-a.s. directional
derivatives for any $v\in V_0(\Lambda)$. 
Let $v\in V_0(\Lambda)$,  any $\omega\in \partial^{c_\zeta}F(\eta)$
must satisfy
\begin{equation*}
  F(({\mathcal V}_t^v)^*\eta)-F(\eta)\le c(({\mathcal V}_t^v)^*\eta+\zeta,\, \omega)- c(\eta+\zeta,\, \omega),
\end{equation*}
for any $t\in \RE$ and $c(\eta+\zeta,\, \omega)<+\infty$. For any
$\beta_{\eta+\zeta,\, \omega}$ realizing $c(\eta+\zeta,\, \omega)$, 
\begin{multline*}
  c(({\mathcal V}_t^v)^*\eta+\zeta,\, \omega) \le \frac{1}{2}\int_{\Lambda \times X}
  \|{\mathcal V}_t^v(x)-y\|^2 \d\beta_{\eta+\zeta,\, \omega} \\
+ \frac{1}{2}\int_{\Lambda^c \times X}
  \|x-y\|^2 \d\beta_{\eta+\zeta,\, \omega}.
\end{multline*}
Hence, 
\begin{equation*}
   F(({\mathcal V}_t^v)^*\eta)-F(\eta)\le \frac{1}{2}\int_{\Lambda \times X}
   (\|{\mathcal V}_t^v(x)-y\|^2 - \|x-y\|^2) \d\beta_{\eta+\zeta,\, \omega}.
\end{equation*}
Divide the two terms of this inequality by $t>0$ and let $t$ go to $0$,
we get 
\begin{equation*}
  \frac{d}{dt}F(({\mathcal V}_t^v)^*\eta)_{|t=0}\le \int_{\Lambda \times X}
  (x-y).v(x)   \d\beta_{\eta+\zeta,\, \omega} (x,y).
\end{equation*}
Applying this inequality to $-v$, we deduce that for any $v\in V_0(\Lambda)$,
\begin{equation*}
  \nabla^\Gamma_v  F(\eta) = \int_{\Lambda \times X}
  (x-y).v(x)   \d\beta_{\eta+\zeta,\, \omega} (x,y).
\end{equation*}
We infer from this relation that for any $\omega \in
\partial^{c_\zeta}F(\eta)$, 
\begin{equation*}
  r^\Lambda(\beta_{\eta+\zeta,\, \omega})=\sum_{x\in\eta}\varepsilon_{(x,\,
    \Id-\nabla^\Gamma _x F(\eta))},
\end{equation*}
 for any
$\beta_{\eta+\zeta,\, \omega}$ realizing $c(\eta+\zeta,\, \omega).$
\end{proof}
The last theorem  means that only a part of any element $\omega$ of
$\partial^{c_\zeta}F(\eta)$ is uniquely determined, namely the part
which will be married to $\eta$ in an optimal coupling between
$\omega$ and $\eta+\zeta$. Nonetheless, when $\zeta=\emptyset$, this
means that $\partial^{c_\emptyset}(\eta)$ is reduced to one point
which is $(\Id -\nabla^\Gamma F)^\Gamma(\eta)=\sum_{x\in\eta}\varepsilon_{x-\nabla^\Gamma _x F(\eta)}$.  

\begin{corollary}
  \label{cor:symmetric_convex_cost}
  Assume that $\hat{\mu}_n$ and $\hat{\nu}_n$ are two absolutely
    continuous, symmetric, probability measures on
  $\tilde{\Lambda}^n$  and that $\T_{c_n^s}(\hat{\mu}_n,\, \hat{\nu}_n)$ is finite. Then there
  exists a unique optimal measure $\rho_n$  for MKP$(\hat{\mu}_n,\, \hat{\nu}_n, \,
  c_n^s)$ and there exists a unique map $\hat{t}_n$ such that $\rho_n=(\Id
  \, \otimes \, \hat{t}_n)^*\hat{\mu}_n$.
\end{corollary}
\begin{proof}
View $\tilde{\Lambda}^n$ as a subset of the Polish space $\Lambda^n$.
Since $\Lambda^n\backslash \tilde{\Lambda}^n$ has a null Lebesgue
measure, we can then view  $\hat{\mu}_n$ and $\hat{\nu}_n$ as absolutely
continuous, symmetric, probability measures on $\Lambda^n$. 
  Since $\Lambda^n$ is Polish, there exists at least one optimal
  measure for MKP$(\hat{\mu}_n,\, \hat{\nu}_n, \,
  c_n^s)$. For any optimal measure $\rho$, there exists a
  $c_n^s$-concave function $f_n$ such that $\supp \rho
  \subset \partial_{c_n^s}f_n$. According to Corollary
  \ref{cor:single-valued}, $\partial_{c_n^s}f_n$ is
  $\hat{\mu}_n$-a.s. single-valued, 
  hence the uniqueness of $\rho_n$ and $t_n$ follows from  Theorem \ref{thm:single_valued_map_implies_uniqueness}.
\end{proof}
We can then state:
\begin{theorem}
\label{thm:transport_fpp}
  Let $\mu$ be a regular probability measure on $\Gamma_\Lambda$ and
  $\nu$ be a probability measure on $\Gamma_X$. The Monge-Kantorovitch
  distance, associated to $c$, between $\mu$ and $\nu$ is finite if
  and only if the following two conditions hold
  \begin{enumerate}[(a) ]
  \item \label{item:4} $\mu(\eta(\Lambda)=n)=\nu(\omega(X)=n)$ for any integer
  $n\ge 0$,
\item  \label{item:5} $\sum_{n\ge 1}\T_{c}(\mu_n,\, \nu_n)^2\mu(\eta(\Lambda)=n)$ is finite.
  \end{enumerate}
Moreover, we have 
\begin{equation}
\label{eq:8}
  \T_c(\mu,\, \nu)^2=\sum_{n\ge 1}\T_{c}(\mu_n,\, \nu_n)^2\mu(\eta(\Lambda)=n),
\end{equation}
and there exists a unique $c$-concave map $F$ such
  that the unique optimal measure $\rho$  is given by
  \begin{equation*}
    \rho=(\Id\  \otimes\  (\Id -\nabla^\Gamma F)^\Gamma)^*\mu,
  \end{equation*}
where 
\begin{equation*}
  (\Id -\nabla^\Gamma F)^\Gamma(\eta)=\sum_{x\in \eta} \varepsilon_{x -\nabla_x^\Gamma F(\eta)}.
\end{equation*}
\end{theorem}
\begin{proof}
  If $\T_c(\mu,\, \nu)$ is finite then according to Theorem
  \ref{thm:egalite_loi_nb_atomes}, condition (\ref{item:4}) is
  satisfied.  Thus, we have
  \begin{align*}
    \T_c(\mu,\, \nu)^2&=\inf_{\gamma\in \Sigma(\mu,\, \nu)}\sum_{n\ge 1}\int_{\eta(\Lambda)=\omega(X)=n}c(\eta,\, \omega)\d\gamma(\eta,\, \omega)\\
& = \inf_{\gamma\in \Sigma(\mu,\, \nu)}\sum_{n\ge 1}\int_{\Gamma_\Lambda^n\times \Gamma_X}c(\eta,\, \omega)\d(\gamma\, |\, \eta(\Lambda)=n)(\eta,\, \omega)\ \mu(\eta(\Lambda)=n)\\
&=\sum_{n\ge 1} \inf_{\gamma_n\in \Sigma(\mu_n,\, \nu_n)}\int_{\Gamma_\Lambda^n\times \Gamma_X}c(\eta,\, \omega)\d\gamma_n(\eta,\, \omega)\ \mu(\eta(\Lambda)=n)\\
&=\sum_{n\ge 1} \T_{c_n^s}(\hat{\mu_n},\, \hat{\nu_n})^2\ \mu(\eta(\Lambda)=n)\\
&=\sum_{n\ge 1} \T_{c}(\mu_n,\, \nu_n)^2\ \mu(\eta(\Lambda)=n),
  \end{align*}
  where $\hat{\mu_n}$ (resp. $\hat{\nu_n}$) is the symmetric
  measure on $\tilde{\Lambda}^n$ corresponding to $\mu_n$
  (resp. $\hat{\mu_n}$). Let $\rho$ an optimal measure whose
  existence is guaranteed because $\Gamma_\Lambda$ and $\Gamma_X$ are
  Polish, we infer from Theorem \ref{thm:compact-case} that there
  exists a $c^\emptyset$-concave function $F$ whose
  $c^\emptyset$-super-gradient is $\mu$-a.s. single valued such that
  $\supp \rho\subset \partial^{c_\emptyset}F$. According to Theorem
  \ref{thm:single_valued_map_implies_uniqueness}, this implies that
  $\rho$ and $T$ are unique and that $\rho=(\Id\  \otimes\ T)^*\mu$. At last, Theorem
  \ref{thm:compact-case} entails that $T=(\Id
  -\nabla^\Gamma F)^\Gamma$.

In the converse direction, since $\mu$ is regular and $\T_{c}(\mu_n,\,
\nu_n)=\T_{c_n^s}(\hat{\mu_n},\, \hat{\nu_n})$ is finite for
any $n\ge 1$, there exists, for any $n\ge 1$,  according to Corollary
\ref{cor:symmetric_convex_cost}, a measure $\hat{\rho}_n$ such that 
\begin{equation*}
  \T_{c_n^s}(\hat{\mu_n},\, \hat{\nu_n})^2=\int_{\Lambda\times X} c_n^s(x,\,y
  )\d\rho_n(x,\, y).
\end{equation*}
Now, we set $$\rho(A)=\sum_{n\ge 1}\rho_n\bigl( A\, \cap
(\Gamma_\Lambda^{(n)}\times \Gamma_X)\bigr)\mu(\eta(\Lambda)=n).$$
Since $\mu(\eta(\Lambda)=n)=\nu(\eta(\Lambda)=n)$ and since
$\hat{\rho}_n$ belongs to $\Sigma(\mu_n, \, \nu_n)$, it is clear that
$\rho$ belongs to $\Sigma(\mu, \, \nu)$. Moreover, we have:
\begin{align}
  \int_{\Gamma_\Lambda\times \Gamma_X}c(\eta,\, \omega)\d \rho(\eta,\,
  \omega) & = \sum_{n\ge 1} \mu(\eta(\Lambda)=n)\int_{\Gamma_\Lambda^{(n)}\times \Gamma_X}c(\eta,\, \omega)\d \rho_n(\eta,\,
  \omega)\nonumber\\
&= \sum_{n\ge 1} \mu(\eta(\Lambda)=n)\int_{\Lambda^n\times
  X^n}c_n^s(x,\, y)\d \hat{\rho}_n(x, \, y)\nonumber\\
&= \sum_{n\ge 1} \mu(\eta(\Lambda)=n) \T_{c_n^s}(\hat{\mu_n},\, \hat{\nu_n})^2,\label{eq:7}
\end{align}
and the last quantity is finite according to the hypothesis. Thus,
$\tau_c(\mu, \, \nu)$ is finite. It remains to prove that $\rho$
constructed above is optimal. For, remind that, as mentioned in the
preliminaries, it is sufficient that $\supp \rho$ be $c$-cyclically monotone.
We infer from the finiteness of $\int c \d \rho$ that for any
$(\eta,\, \omega)$ in $\supp \rho$, $\eta(\Lambda)=\omega(X)$. For $m$
any integer, let
$\bigl((\eta_i,\, \omega_i),\, i=1, \cdots, m\bigr)$  be a family of
elements of $\supp \rho$. Set $I_n=\{i\in 1,\cdots,\, m,\
\eta_i(\Lambda)=n\}$, we can then write
\begin{equation*}
  \sum_{i=1}^m c(\eta_i,\, \omega_i) = \sum_{n=1}^{+\infty}
  \sum_{i\in I_n}c(\eta_i,\, \omega_i).
\end{equation*}
Let $\sigma\in \SS_m$, if for some $n$, $\sigma I_n$ differs from
$I_n$ then $\sum_{i\in I_n}c(\eta_i,\, \omega_i)$ is infinite and it
is clear that 
$$ \sum_{i=1}^m c(\eta_i,\, \omega_i)\le  \sum_{i=1}^m c(\eta_i,\,
\omega_{\sigma(i)}).$$
Thus, we now assume that for any $n\ge 1$, $\sigma I_n=I_n$, i.e.,
$\omega_i(\Lambda)=\omega_{\sigma(i)}(\Lambda)$ for any $i=1,
\cdots,\, m$. Since for any $n\ge 1$, $\hat{\rho}_n$ is
$c_n^s$-cyclically monotone, so does $\rho_n$. Moreover, $\supp \rho
\cap I_n=\rho_n$, thus for any $n\ge 1$,
\begin{equation*}
   \sum_{i\in I_n}c(\eta_i,\, \omega_i)\le  \sum_{i\in I_n}c(\eta_i,\, \omega_{\sigma(i)}).
\end{equation*}
By summation, we infer that $\sum_{i=1}^m c(\eta_i,\, \omega_i)\le
\sum_{i=1}^m c(\eta_i,\, \omega_{\sigma(i)})$ for any $\sigma \in
\SS_m$. This amounts to say that $\supp \rho$ is $c$-cyclically
monotone, hence that $\rho$ is an optimal measure (unique according to
the first part of the proof) for  MKP$(\mu,\, \nu, \,  c)$. We deduce 
from (\ref{eq:7}) that (\ref{eq:8}) holds true. 
\end{proof}

\subsection{Example : Wasserstein distance with respect to a Poisson process}
\label{sec:example-:-wass}
Let $\sigma$ be a diffuse (by which we mean absolutely continuous with
respect to the Lebesgue measure) Radon measure on $X,$ the Poisson
measure on $\Gamma_X$ with intensity $\sigma$, denoted by $\mu_\sigma,$ is the unique
probability measure on $(  \Gamma_X,{\mathcal B}(\Gamma_X))$ such that
\begin{equation}
  \label{eq:1}
  E[\exp(\int f \d \eta)]=\exp\left(\int_{X} (e^{f(x)}-1)\d \sigma(x)\right),
\end{equation}
for all $f\in{\mathcal C}_0.$ It is well known  that $\mu_\sigma$
satisfies Hypothesis \ref{hyp:hyp} and $\mu_\sigma$ is regular since
$\mu_n=\sigma^{\otimes n}$,  thus
the previous results apply.
Let $\sigma_1$ and $\sigma_2$ two diffuse probability measures on
 $X$ with finite Wasserstein distance with respect to the euclidean
 cost on $X=\R^k$:
 \begin{equation*}
   \T_e(\sigma_1,\, \sigma_2)^2=\inf\limits_{\gamma \in
     \Sigma(\sigma_1,\, \sigma_2)}\frac{1}{2}\int_{X\times X}
   \|x-y\|^2\d \gamma(x,\, y)<+\infty.
 \end{equation*}
 We  denote by $t$ the transport map from $\sigma_1$ to $\sigma_2$ and $\phi$ its potential, i.e., the
 convex map from $X\to\R$ such that $\nabla \phi=t.$ By $\nabla,$ we
 mean here the usual gradient in $X.$ For $W$ and $Y$ two spaces
 and $f:W\to \R$ and $g:Y\to \R,$ we denote by $f\oplus g$ the map
 defined on $W\times Y$ by $(f\oplus g)(x,y)=f(x)+g(y).$
 \begin{lemma}
   The map $t^{\otimes (n)}:\, (x_1,\ldots,x_n)\in X^n\mapsto
   (t(x_1),\ldots,t(x_n))$ is the transport map from $\sigma_1^{\otimes
     n}$ to $\sigma_2^{\otimes n}.$ Moreover,
   \begin{equation*}
%\label{eq:11}
     \T_e(\sigma_1^{\otimes n},\, \sigma_2^{\otimes n})^2=n  \T_e(\sigma_1,\, \sigma_2)^2.
   \end{equation*}
 \end{lemma}
 \begin{proof}
   It is immediate that $t^{\otimes(n)}\equiv \nabla (\oplus_{i=1}^n\phi)$ and
   that $\oplus_{i=1}^n\phi$ is convex, thus $t^{\otimes(n)}$ is cyclically
   monotone (with respect to the squared euclidean cost on $X^n$). Moreover,
   $(t^{\otimes(n)})^*(\nu_1^{\otimes n})=\nu_2^{\otimes n},$ hence $t^{\otimes(n)}$ is the
   optimal transport map between $\nu_1^{\otimes n}$ and
   $\nu_2^{\otimes n}$. Then,
   \begin{align*}
      \T_e(\sigma_1^{\otimes n},\, \sigma_2^{\otimes n})^2
      &= \int_{X^n}\frac{1}{2}\|x-t^{\otimes(n)}(x)\|^2 \d
      \sigma_1^{\otimes n}(x)\\
& =\sum_{j=1}^n \int_X  \frac{1}{2}\|x_j-t^{\otimes(n)}(x_j)\|^2 \d
      \sigma_1(x_j)\\
 &= n\T_e(\sigma_1,\, \sigma_2)^2.
   \end{align*}
The proof is thus complete.
 \end{proof}
It then follows from Theorem \ref{thm:transport_fpp} that: 
 \begin{theorem}\label{thm:transport_poisson}
For $\sigma_1$ and $\sigma_2$ two diffuse probability measures on
 $X$, if $\T_e(\sigma_1, \, \sigma_2)<+\infty$ then
 $\T_c(\mu_{\sigma_1},\, \mu_{\sigma_2})$ is finite. If
 $t=\nabla \phi$ are respectively the transport map and its associated
 potential for MKP$(\sigma_1,\, \sigma_2,\, c_e)$ then 
 \begin{equation*}
   T=\sum_{n\ge 1} t^{\otimes(n)}\car_{\Gamma_X^{(n)}} \text{ and }
   \Phi =\sum_{n\ge 1} (\oplus_{i=1}^n\phi)^\Gamma\car_{\Gamma_X^{(n)}} ,
 \end{equation*}
are respectively the transport map and the associated potential for
the Monge-Kantorovitch problem MKP$(\mu_{\sigma_1},\, \mu_{\sigma_2},\, c)$.
Moreover,
\begin{equation}
\label{eq:14}
  \T_c(\mu_{\sigma_1},\, \mu_{\sigma_2})=\T_e(\sigma_1, \, \sigma_2).
\end{equation}
 \end{theorem}
 \begin{rem}
   For finite point processes, it is possible to define a cost between
   configurations by 
   \begin{equation*}
     c_b(\eta, \, \omega)=\frac{1}{\eta(X)}c(\eta,\, \omega). 
   \end{equation*}
We would then have 
\begin{equation*}
  \T_{c_b}(\mu_{\sigma_1},\, \mu_{\sigma_2})^2=(1-e^{-1})\T_e(\sigma_1, \, \sigma_2)^2.
\end{equation*}
This  distance $\T_{c_b}$  appears in  papers of Barbour et al. \cite{MR1190904,MR2001a:60058}.
 \end{rem}
A Cox process is a doubly-stochastic Poisson process: $\sigma$ is now
a random variable in the set of diffuse Radon measures on $X$ and
conditionally to $\sigma$, the point process is a Poisson process of
intensity $\sigma$. By conditioning with respect to the intensities,
the proof given above yields to the following theorem.
 \begin{theorem}
   If $\mu$ and $\nu$ are two Cox processes of random intensities
   $\sigma_1$ and $\sigma_2$ respectively, such that
   $\esp{\T_e(\sigma_1,\sigma_2)}$ is finite. Then,
   \begin{equation*}
     \T_c(\mu_{\sigma_1},\mu_{\sigma_2})=\esp{\T_e(\sigma_1,\sigma_2)}.
   \end{equation*}
 \end{theorem}

%%%%%%%%%%%%%%%%%%%%%%%%%%%%%%%%%%%%%%%%%%%%%%%%%%%%%%%%%%%%%%%%%%%%%%%%%%%%%%%%%%%%%%%%%%%%%%%%
%%%%%%%%%%%%%%%%%%%%%%%%%%%%%%%%%%%%%%%%%%%%%%%%%%%%%%%%%%%%%%%%%%%%%%%%%%%%%%%%%%%%%%%%%%%%%%%%
\section{Locally finite point processes}
\label{sec:locally-finite-point}We now only assume that $\mu$ is the
law of a locally finite point
process : $\mu(\eta(\Lambda)<+\infty)=1$ for all compact sets
$\Lambda$ but $\mu(\eta(X)=+\infty)=1$. We can no longer work on the graded
space $\cup_{n\ge 1}\Gamma_{X}^{(n)}$ since it is
$\mu$-negligible. We are in fact reminded the case of the Wiener
space. There is thus no big surprise that we can follow closely the
beautiful method of \cite{MR2036490}.
\begin{lemma}
\label{lem:optimalite_prob_conditionne}
  Let $\mu$ and $\nu$ be two probability measures on $\Gamma_X$ such
  that $\T_c(\mu,\, \nu)$ is finite. Let $\gamma$ be one optimal measure
  and $\Lambda$ be a compact set of $X$. Consider the disintegration
  of $\gamma$ along the projection $\pi_1^\Lc$, i.e.,
  \begin{equation*}
    \gamma(.)=\int_{\Gamma_\Lc} \gamma(.\,  | \, \pi_1^\Lc(\eta,\, \omega)=\eta_\Lc)\d \mu_\Lc(\eta_\Lc),
  \end{equation*}
where $\mu_\Lc$ is the image measure of $\mu$ by $\pi^\Lc$. Denote by
$\gamma(.\, |\, \eta_\Lc)$ the  regular version of the conditional
probability $\gamma(.\,  | \, \pi_1^\Lc(\eta,\, \omega)=\eta_\Lc)$. Then,
$\mu_\Lc$-a.s., $\gamma(.\, |\, \eta_\Lc)$ is an optimal measure for
MKP$(p_1\gamma(.\, |\, \eta_\Lc), \ p_2\gamma(.\, |\, \eta_\Lc), c_{\eta_\Lc})$.
\end{lemma}
\begin{rem}
  If we denote by $(N,\, M)$ a couple of random variables whose
  distribution is $\gamma$ and if $N_\Lambda(\eta,\omega):=N(\eta)\cap
  \Lambda$, then the previous lemma stands that conditionally to
  $(N_{\Lambda^c}=\eta_{\Lambda^c})$, the law of
  $(\eta_{\Lambda^c}+N_\Lambda,\, M)$ is optimal for MKP($\gamma_{
    \eta_{\Lambda^c}+N_\Lambda\,|\, N_{\Lambda^c}=\eta_{\Lambda^c}},\
  \gamma_{M\,|\, N_{\Lambda^c}=\eta_{\Lambda^c}},\,
  c_{\eta_\Lambda^c}$).
Note that within this setting, since the law of $N$ is $\mu$, it is clear that
\begin{equation*}
 \gamma_{N \,|\,  N_{\Lambda^c}=\eta_{\Lambda^c}}=\mu_{N \, | \,
  N_{\Lambda^c}=\eta_{\Lambda^c}} \text{ i.e., } p_1\gamma(.\, |\, \eta_\Lc)=\mu(.\, | \, \eta_\Lc).
\end{equation*}
  \end{rem}
  \begin{proof}[Proof of Lemma
    \protect\ref{lem:optimalite_prob_conditionne}]
According to the  definition of an optimal measure,
\begin{align*}
  J_c(\gamma)&=\int c(\eta,\, \omega)\d\gamma(\eta,\, \omega)\\
&= \int c(\pi^\Lambda\eta +  \pi^\Lc\eta,\, \omega)\d\gamma(\eta,\, \omega)\\
&=\int_{\Gamma_\Lc}\d\mu_\Lc(\eta_\Lc)\ \int_{\Gamma_X\times \Gamma_X}
c(\eta_\Lc + \pi^\Lambda \eta,\, \omega)\d\gamma(\eta,\, \omega\, | \,
\eta_\Lc)\\
&= \int_{\Gamma_\Lc}\d\mu_\Lc(\eta_\Lc)\ \int_{\Gamma_\Lambda\times \Gamma_X}
c_{\eta_\Lc} ( \eta,\, \omega) \d(\pi^\Lambda\otimes\Id)\gamma(\eta,\,
\omega\, | \,
\eta_\Lc)\\
&=J_{c_{\eta_\Lc}}((\pi^\Lambda\otimes\Id)\gamma(. \, | \,
\eta_\Lc)).
\end{align*}
Now, note that $(\pi^\Lambda\otimes\Id)\gamma(. \, | \,
\eta_\Lc)$ has marginals $\pi^\Lambda\mu(.\, | \, \eta_\Lc)$ and
$p_2\gamma(. \, | \, \eta_\Lc)$ which are probability measures on
$\Gamma_\Lambda$ and $\Gamma_X$ respectively. Let ${\mathcal
  M}_1(\Gamma_\Lambda \times \Gamma_X)$ be the set of probability
measures on $\Gamma_\Lambda\times \Gamma_X$. Define the sets $B$ and
$C$ as
\begin{gather*}
  B=\left\{(\eta_\Lc,\, \theta):\ \theta\in \Sigma(\pi^\Lambda\mu(.\, | \,
  \eta_\Lc), \ p_2\gamma(. \, | \, \eta_\Lc)\right\}\\
C=\left\{(\eta_\Lc,\, \theta):\ J_{c_{\eta_\Lc}}(\theta) <
J_{c_{\eta_\Lc}}(\gamma(. \, | \, \eta_\Lc))\right\}.
\end{gather*}
Let $K$ be the projection on $\Gamma_\Lc$ of $C$. Since $B$ and
$C$ are Borel, $K$ is a Souslin set, hence $\mu_\Lc$-measurable. Thus
there exists a measurable map $\Theta$ from $K$ to  ${\mathcal
  M}_1(\Gamma_\Lambda \times \Gamma_X)$ such that $(\eta_\Lc,\,
\Theta(\eta_\Lc))$ belongs to $C$, for $\mu_\Lc$-almost-all
$\eta_\Lc$. Define a measure $\theta$ as:
\begin{equation*}
  \theta = \int_K
  \Theta(\eta_\Lc)\d\mu_\Lc(\eta_\Lc)+\int_{K^c}(\pi^\Lambda\otimes \Id)\gamma(.\, | \, \eta_\Lc)\d\mu_\Lc(\eta_\Lc).
\end{equation*}
If $\mu_\Lc(K)>0$ then
\begin{multline*}
  J_c(\theta)
=\int_{K^c} J_{c_{\eta_\Lc}}((\pi^\Lambda\otimes\Id)\gamma(. \, | \,
\eta_\Lc))\d\mu_\Lc(\eta_\Lc) \\
\shoveright{+ \int_{K}
J_{c_{\eta_\Lc}}(\Theta(\eta_\Lc))\d\mu_\Lc(\eta_\Lc)}\\
\shoveleft{<  \int_{K^c}
J_{c_{\eta_\Lc}}((\pi^\Lambda\otimes\Id)\gamma(. \, | \,
\eta_\Lc))\d\mu_\Lc(\eta_\Lc)}\\
\shoveright{ +\int_{K}
J_{c_{\eta_\Lc}}((\pi^\Lambda\otimes\Id)\gamma(. \, | \,
\eta_\Lc))\d\mu_\Lc(\eta_\Lc)}\\
=J_c(\gamma),    
\end{multline*}
which is a contradiction to the optimality of $\gamma$.
  \end{proof}
  \begin{theorem}
    \label{thm:lemme2}
Assume that the hypothesis of Lemma
\ref{lem:optimalite_prob_conditionne} holds and assume that $\mu$ is
 regular. Let $\Lambda$ be any
compact subset of $X$. Then, there exists a measurable map $T_\Lambda$ from
$\Gamma_X$ to itself such that 
\begin{equation*}
  \gamma\left(
(\eta,\, \omega): \, r^\Lambda_2(\{\beta_{\eta,\,
  \omega}\})=T_\Lambda(\pi^\Lambda\eta, \, \pi^\Lc\eta)\,\right)=1
\end{equation*}
  \end{theorem}
  \begin{proof}
    Fix $\eta_\Lc \in \Gamma_\Lc$ and define $C_{\eta_\Lc}$ as the
    support of $\gamma( . \, | \, \eta_\Lc)$. Consider the 
%following
%    diagram
%    \begin{equation*}
%      \begin{array}{ccccc}
%      C_{\eta_\Lc}  & \longrightarrow & K_{\eta_\Lc} & \longrightarrow
%      & K_{\eta_\Lc,\, \eta}\\
%(\eta , \, \omega) & \longmapsto & (\eta+\eta_\Lc, \,
%r^\Lambda_2(\{\beta_{\eta+\eta_\Lc,\ \omega}\})&&\\
%&&(\eta, \, \zeta) & \longmapsto & \zeta
%      \end{array}
%    \end{equation*}
%This means that 
two sets:
\begin{equation*}
  K_{\eta_\Lc}=\left\{(\eta, \,r^\Lambda_2(\{\beta_{\eta+\eta_\Lc,\
      \omega}\}))\in \Gamma_\Lambda\times\Gamma_X:\ (\eta+\eta_\Lc, \,
 \omega\})\in C_{\eta_\Lc}
\right\}
\end{equation*}
and
\begin{equation*}
  K_{\eta_\Lc,\, \eta}=\left\{ \omega\in \Gamma_X, \, (\eta,\omega)\in
     K_{\eta_\Lc}\right\}.
\end{equation*}
We know from Theorem \ref{thm:compact-case} that for $\mu_\Lc$-almost
all $\eta_\Lc$, $K_{\eta_\Lc,\,
  \eta}$ is reduced to one point for $\mu_\Lambda(. \, | \,
\eta_\Lc)$-almost-all $\eta$. Let 
\begin{equation*}
  N=\{(\eta,\, \eta_\Lc)\in \Gamma_\Lambda\times \Gamma_\Lc:\
  \text{ Card }(K_{\eta,\, \eta_\Lc})>1\},
\end{equation*}
$N$ is a Souslin set, hence it is universally measurable. Let $\sigma$
be the measure defined as the image of $\mu$ under the projection
$\eta\mapsto (\pi^\Lambda\eta,\, \pi^\Lc\eta)$. We then have
\begin{equation*}
  \sigma(N)=\int_{\Gamma_\Lc}\d\mu_\Lc(\eta_\Lc)\
  \int_{\Gamma_\Lambda}\car_N(\eta,\, \eta_\Lc) \mu(\d\eta\, | \, \eta_\Lc)=0.
\end{equation*}
Hence, $\mu$ and $\gamma$ almost-surely, $ K_{\eta_\Lc,\, \eta}$ is
reduced to a single-point and we define $T_\Lambda$ as the map which
sends $(\eta,\, \eta_\Lc)$ to this point. It is automatically
measurable by the selection theorem.
  \end{proof}
  \begin{theorem}
Assume that the hypothesis of Lemma
\ref{lem:optimalite_prob_conditionne} holds and assume that $\mu$ is
 regular. Let $(\Lambda_n, \, n\ge 1)$ be an increasing sequence of
 compact sets such that $\cup_{n\ge 1}\Lambda_n=X$. Then, there exists
 a unique optimal measure $\gamma$ and a unique map $T$ such that 
 \begin{equation*}
   \gamma=(\Id \otimes T)^* \mu.
 \end{equation*}
  \end{theorem}
  \begin{proof}
Let $\gamma$ be an optimal measure for MKP($\mu, \, \nu, \, c$).
According to Theorem \ref{thm:lemme2}, we know that 
\begin{equation*}
  r^{\Lambda_n}_2(\{\beta_{\eta,\,
  \omega}\})=T_{\Lambda_n}(\pi^{\Lambda_n}\eta, \, \pi^\Lnc\eta)
\end{equation*}
$\gamma$-a.s. for all integers $n$.
Let $B$ be a bounded subset of $X$, we clearly have
\begin{equation*}
   r^{\Lambda_n}_2(\beta_{\eta,\,
  \omega})(B)\le \omega(B)<\infty,
\end{equation*}
for any $\beta_{\eta,\,
  \omega}$ realizing $c(\eta, \, \omega)$. Thus, for $\gamma$-almost
all $(\eta, \, \omega)$, the family $( r^{\Lambda_n}_2(\{\beta_{\eta,\,
  \omega}\}), \, n\ge 1)$ is tight in $\Gamma_X$ (see \cite{kallenberg83}). Hence, up to the
extraction of a subsequence, one can assume that $ r^{\Lambda_n}_2(\{\beta_{\eta,\,
  \omega}\})$ converges to  $\omega$. On
the other hand, $\pi^{\Lambda_n}\eta$ converges to  $\eta$ and
$\pi^\Lnc\eta$ converges to $\emptyset$ as $n$ goes to
infinity. Define $T$ by $T(\eta)=\lim_{n\to
  \infty}T_{\Lambda_n}(\pi^{\Lambda_n}\eta, \, \pi^\Lnc\eta)$, we
clearly have $\omega=T(\eta)$, $\gamma$-a.s.
The conclusion follows by Theorem \ref{thm:single_valued_map_implies_uniqueness}.
  \end{proof}
We didn't manage to find any sufficient condition which would ensure
the finiteness of the Wasserstein distance between two locally finite
point processes. However, we do know that there exists some relevant
cases. Consider, for instance, we are given a Poisson process of non-finite
intensity $\sigma_1$ and a map $h$ from $X$ to itself such that $\int \|h\|^2 \d \sigma_1$ is finite.
Then, 
\begin{align*}
  \T_{c}(\mu_{\sigma_1},\, (\Id +
  h^\Gamma)^*\mu_{\sigma_1})&\le\frac{1}{2}  \esp{\sum_{x\in \eta}\|
    x-(\Id+h)(x)\|^2\d \eta(x)}\\
&=\frac{1}{2}  \int \|h\|^2 \d \sigma_1 < \infty.
\end{align*}
Note that $ (\Id +
  h^\Gamma)^*\mu_{\sigma_1}$ is a Poisson process of intensity $(\Id +h)^*\sigma_1$.
%%%%%%%%%%%%%%%%%%%%%%%%%%%%%%%%%%%%%%%%%%%%%%%%%%%%%%%%%%%%%%%%%%%%%%%%%%%%%%%%%%%%%%%%%%%%%%%%
%\bibliography{references}

\begin{thebibliography}{AKR98}

\bibitem[AKR98]{MR99d:58179}
S.~Albeverio, Yu.~G. Kondratiev, and M.~Rockner, \emph{Analysis and geometry on
  configuration spaces}, J. Funct. Anal. \textbf{154} (1998), no.~2, 444--500.
  \MR{99}

\bibitem[BB92]{MR1190904}
A.~D. Barbour and T.~C. Brown, \emph{Stein's method and point process
  approximation}, Stochastic Process. Appl. \textbf{43} (1992), no.~1, 9--31.
  \MR{MR1190904 (93k:60120)}

\bibitem[BC01]{MR2002k:60029}
A.~D. Barbour and O.~Chryssaphinou, \emph{Compound {P}oisson approximation: a
  user's guide}, Ann. Appl. Probab. \textbf{11} (2001), no.~3, 964--1002.
  \MR{2002}

\bibitem[BHJ92]{MR93g:60043}
A.~D. Barbour, L.~Holst, and S.~Janson, \emph{Poisson approximation}, Oxford
  Studies in Probability, vol.~2, The Clarendon Press Oxford University Press,
  1992, Oxford Science Publications. \MR{93}

\bibitem[BM02]{MR1920275}
A.~D. Barbour and M.~M{aa}nsson, \emph{Compound {P}oisson process
  approximation}, Ann. Probab. \textbf{30} (2002), no.~3, 1492--1537. \MR{1}

\bibitem[BX00]{MR2001e:60040}
A.~D. Barbour and A.~Xia, \emph{Estimating {S}tein's constants for compound
  {P}oisson approximation}, Bernoulli \textbf{6} (2000), no.~4, 581--590.
  \MR{2001}

\bibitem[DVJ03]{MR1950431}
D.~J. Daley and D.~Vere-Jones, \emph{An introduction to the theory of point
  processes. {V}ol. {I}}, second ed., Probability and its Applications (New
  York), Springer-Verlag, New York, 2003, Elementary theory and methods.
  \MR{MR1950431 (2004c:60001)}

\bibitem[F{\"U}04]{MR2036490}
D.~Feyel and A.~S. {\"U}st{\"u}nel, \emph{Monge-{K}antorovitch measure
  transportation and {M}onge-{A}mp\`ere equation on {W}iener space}, Probab.
  Theory Related Fields \textbf{128} (2004), no.~3, 347--385. \MR{MR2036490
  (2004m:60121)}

\bibitem[Kal83]{kallenberg83}
O.~Kallenberg, \emph{Random measures}, 3rd ed., Academic Press, 1983.

\bibitem[Lev99]{MR1699061}
V.~Levin, \emph{Abstract cyclical monotonicity and {M}onge solutions for the
  general {M}onge-{K}antorovich problem}, Set-Valued Anal. \textbf{7} (1999),
  no.~1, 7--32. \MR{MR1699061 (2000j:90075)}

\bibitem[RR98a]{MR99k:28006}
S.~T. Rachev and L.~R{\"u}schendorf, \emph{Mass transportation problems. {V}ol.
  {I}}, Probability and its Applications (New York), Springer-Verlag, New York,
  1998, Theory. \MR{99}

\bibitem[RR98b]{MR99k:28007}
\bysame, \emph{Mass transportation problems. {V}ol. {II}}, Probability and its
  Applications (New York), Springer-Verlag, New York, 1998, Applications.
  \MR{99}

\bibitem[RS99]{MR1730565}
M.~R{\"o}ckner and A.~Schied, \emph{Rademacher's theorem on configuration
  spaces and applications}, J. Funct. Anal. \textbf{169} (1999), no.~2,
  325--356. \MR{MR1730565 (2001b:58058)}

\bibitem[R{\"u}s96]{MR1395577}
L.~R{\"u}schendorf, \emph{On {$c$}-optimal random variables}, Statist. Probab.
  Lett. \textbf{27} (1996), no.~3, 267--270. \MR{MR1395577 (97h:62051)}

\bibitem[Tho00]{MR1741181}
H.~Thorisson, \emph{Coupling, stationarity, and regeneration}, Probability and
  its Applications (New York), Springer-Verlag, New York, 2000, pp.~xiv+517.

\bibitem[Vil03]{MR1964483}
C.~Villani, \emph{Topics in optimal transportation}, Graduate Studies in
  Mathematics, vol.~58, American Mathematical Society, Providence, RI, 2003.
  \MR{MR1964483 (2004e:90003)}

\bibitem[Xia00]{MR2001a:60058}
A.~Xia, \emph{Poisson approximation, compensators and coupling}, Stochastic
  Anal. Appl. \textbf{18} (2000), no.~1, 159--177. \MR{2001}

\end{thebibliography}
%\bibliographystyle{amsalpha}
\providecommand{\bysame}{\leavevmode\hbox to3em{\hrulefill}\thinspace}
\providecommand{\MR}{\relax\ifhmode\unskip\space\fi MR }
% \MRhref is called by the amsart/book/proc definition of \MR.
\providecommand{\MRhref}[2]{%
  \href{http://www.ams.org/mathscinet-getitem?mr=#1}{#2}
}
\providecommand{\href}[2]{#2}

\end{document}